# M.V. KHARINOV
# PRODUCT OF THREE OCTONIONS

*Khariniv M.V.* **Product of three octonions.**

**Abstract**. This paper is devoted to octonions that are the eight-dimensional hypercomplex numbers characterized by multiplicative non-associativity. The decomposition of the product of three octonions with the conjugated central factor into the sum of mutually orthogonal anticommutator, commutator and associator, is introduced in an obvious way by commuting of factors and alternating the multiplication order. The commutator is regarded as a generalization of the cross product to the case of three arguments both for quaternions and for octonions. It is verified that the resulting additive decomposition is equivalent to the known solution derived and presented by S. Okubo in a cumbersome form.
**Keywords:** quaternions, octonions, triple cross product, additive decomposition, mutually orthogonal terms.

**1. Introduction.** The famous generalization of real numbers by W.R. Hamilton surpassed by time the introduction of conventional cross product in mathematics and the discovery of four-dimensional space-time in physics [1]. Nowadays many applications of hypercomplex numbers are spawned not only in the field of classical physics and mathematics, but also in variety branches of modern science and technology, where the utilization of quaternions is constantly expanding. Nevertheless, according to opinions of many experts, the potential of four-dimensional quaternions and especially capabilities of their eight-dimensional generalization (octonions) are far from been exhausted.

Unlike the quaternions, the treatment of non-associative octonions, although occurs [2], but comparatively much less frequently, since it is necessary to deal with non-associative vector multiplication. A possible way to solve the problem is to generalize the basic method of additive decomposition of the product of hypercomplex numbers into mutually orthogonal symmetric and antisymmetric components, which is the topic of the paper.

The paper is intended for a reader familiar with the hypercomplex numbers (normalized algebras with multiplicative unit). Otherwise, the general information on hypercomplex numbers can be obtained from Wikipedia. All the necessary information in an accessible form is set forth in the popular book [3]. Many useful rules for working with hypercomplex numbers are given in [4]. In this paper we list only the most relevant formulas. Unnumbered formulas are given mostly to explain the meaning of notations.

**2. Elementary information.** By $i_0$ we denote the multiplicative unit that commutes with any hypercomplex number and, when multiplied, leaves it unchanged:
$$i_0 u = u i_0 = u.$$

Let $\bar{u}$ denotes the hypercomplex number which is conjugate to the number $u$, and is related to $u$ by the formula [4]:
$$\bar{u} = 2(u, i_0)i_0 - u,$$
wherein $(u, i_0)$ is the inner product of vector $u$ and vector $i_0$.

In four algebras of hypercomplex numbers (real numbers, complex numbers, quaternions and octonions), the square of the length $(u, u)$ of the vector $u$ is introduced as the product of the vector $u$ and the conjugate vector $\bar{u}$, and the inner product of the vectors $u_1$ and $u_2$ coincides with the half-sum of the vector $u_1 \bar{u}_2$ and the conjugate vector $u_2 \bar{u}_1$:

$$u\bar{u} \equiv \bar{u}u = (u,u)i_0 \quad \Leftrightarrow \quad \frac{u_1 \bar{u}_2 + u_2 \bar{u}_1}{2} \equiv \frac{\bar{u}_1 u_2 + \bar{u}_2 u_1}{2} = (u_1, u_2)i_0.$$

From the last relations it is easy to establish the useful equality:
$$u_1 \bar{u}_2 u_1 = 2(u_1, u_2)u_1 - (u_1, u_1)u_2, \tag{1}$$
wherein the product $u_1 \bar{u}_2 u_1$ is written without brackets, since it does not depend on the order of multiplication of hypercomplex factors:
$$u_1 \bar{u}_2 u_1 = (u_1 \bar{u}_2)u_1 = u_1(\bar{u}_2 u_1).$$

In order to reliably relate to multiplication of quaternions and octonions, it is useful to keep in mind the following elementary formulas.

So, the next formula expresses that conjugation changes the order of factors:
$$\overline{u_1 u_2} = \bar{u}_2 \bar{u}_1.$$

And the following formula for the inner product $(u_1 u_2, i_0)$ vector $u_1 u_2$ and vector $i_0$ expresses the transfer rule of the factor with simultaneous conjugation, from one to another part of the inner product:
$$(u_1 u_2, i_0) \equiv (u_2 u_1, i_0) = (u_1, \bar{u}_2).$$

Finally, the formula:
$$((u_1 u_2)u_3, i_0) = (u_1(u_2 u_3), i_0) \equiv (u_1 u_2 u_3, i_0)$$
states that the order of multiplying of three hypercomplex factors $u_1$, $u_2$ and $u_3$ does not affect the inner product of $(u_1 u_2)u_3$ and $i_0$. So, the product $u_1 u_2 u_3$ written without brackets implies choosing either of the two alternative ways of their arrangement, just as in left part of (1).

## 3. Additive decomposition of the product of two hypercomplex numbers.

The product $u_1 u_2$ of two hypercomplex numbers $u_1$ and $u_2$ usually is decomposed into the sum of two mutually orthogonal terms $\{u_1, u_2\}$ and $[u_1, u_2]$:

$$u_1 u_2 = \frac{u_1 u_2 + u_2 u_1}{2} + \frac{u_1 u_2 - u_2 u_1}{2} = \{u_1, u_2\} + [u_1, u_2], \qquad (2)$$

wherein $(\{u_1, u_2\}, [u_1, u_2]) = 0$, the anticommutator $\frac{u_1 u_2 + u_2 u_1}{2}$, denoted as $\{u_1, u_2\}$, is a linear combination of the arguments $u_1$, $u_2$, and multiplicative unit $i_0$:

$$\{u_1, u_2\} = \frac{u_1 u_2 + u_2 u_1}{2} = (u_1, i_0) u_2 + (u_2, i_0) u_1 - (u_1, u_2) i_0, \qquad (3)$$

and the commutator $\frac{u_1 u_2 - u_2 u_1}{2}$, which is denoted as $[u_1, u_2]$, is a cross product of a pair of hypercomplex numbers $u_1$ and $u_2$.

Taking into account (3) the product of two hypercomplex numbers $u_1 u_2$, is expressed in expanded form by the formula [5]:

$$u_1 u_2 = (u_1, i_0) u_2 + (u_2, i_0) u_1 - (u_1, u_2) i_0 + [u_1, u_2]. \qquad (4)$$

Cross product $[u_1, u_2]$ of hypercomplex numbers $u_1$ and $u_2$ is orthogonal to multiplicative unit $i_0$ and vanishes to zero when the multiplicative unit $i_0$ is used as one or another argument:

$$([u_1, u_2], i_0) = 0,$$

$$[i_0, u_2] = [u_1, i_0] = 0.$$

In other respects, the cross product $[u_1, u_2]$ preserves the properties of the conventional cross product, which is introduced in three-dimensional space using intuitively perceived "right hand rule". It means that, if we agree to denote by the prime an annihilating of the real component of hypercomplex number $u$:

$$u' = u - (u, i_0) i_0,$$

then the cross product $[u_1, u_2]$ of the vectors $u_1$ and $u_2$ coincides with the conventional cross product $[u'_1, u'_2]$ of the vectors $u'_1$ and $u'_2$ of the three-dimensional subspace that is orthogonal to multiplicative unit $i_0$. In particu-lar, the square of the length of the cross product is expressed by the symmetric determinant from the inner products of said three-dimensional vectors:

$$([u_1,u_2],[u_1,u_2]) = ([u'_1,u'_2],[u'_1,u'_2]) = \det\begin{bmatrix}(u'_1,u'_1) & (u'_1,u'_2)\\(u'_1,u'_2) & (u'_2,u')_2\end{bmatrix}^2.$$

wherein $\det\begin{bmatrix}(u'_1,u'_1) & (u'_1,u'_2)\\(u'_1,u'_2) & (u'_2,u'_2)\end{bmatrix} = (u'_1,u'_1)(u'_2,u'_2) - (u'_1,u'_1)^2$.

At the same time, square of the length of the anticommutator $\{u_1,u_2\}$ is expressed as:

$$(\{u_1,u_2\},\{u_1,u_2\}) = (u_1,u_1)(u_2,u_2) - \det\begin{bmatrix}(u'_1,u'_1) & (u'_1,u'_2)\\(u'_1,u'_2) & (u'_2,u'_2)\end{bmatrix}.$$

So the axiomatic property of *normalized* algebras:
$(u_1 u_2, u_1 u_2) \equiv (\{u_1,u_2\},\{u_1,u_2\}) + ([u_1,u_2],[u_1,u_2]) = (u_1,u_1)(u_2,u_2)$
is fulfilled.

**4. Additive decomposition of the product of three hypercomplex numbers.** The product of three octonions obviously can be written in the form:

$$(u_1\bar{u}_2)u_3 = \frac{(u_1\bar{u}_2)u_3 + (u_3\bar{u}_2)u_1}{2} + \frac{(u_1\bar{u}_2)u_3 - u_3(\bar{u}_2 u_1)}{2} + \frac{u_3(\bar{u}_2 u_1) - (u_3\bar{u}_2)u_1}{2}. \quad (5)$$

NOTE: the conjugation of central factor provides mutual orthogonality of summands and is necessary to obtain laconic relations. As any number is conjugate to conjugated one this doesn't reduce the generality of reasoning.

In this case the triple product $(u_1\bar{u}_2)u_3$ is decomposed into the sum of *mutually orthogonal* terms:

$$(u_1\bar{u}_2)u_3 = \{u_1,u_2,u_3\} + [u_1,u_2,u_3] + \langle u_1,u_2,u_3\rangle, \quad (6)$$

wherein
$(\{u_1,u_2,u_3\},[u_1,u_2,u_3]) = (\{u_1,u_2,u_3\},\langle u_1,u_2,u_3\rangle) = (\{u_1,u_2,u_3\},\langle u_1,u_2,u_3\rangle) = 0$.

The mutually orthogonal summands $\{u_1,u_2,u_3\}$, $[u_1,u_2,u_3]$ and $\langle u_1,u_2,u_3\rangle$ are determined by the relations (7)-(9) with an alternative order of multiplication of the factors:

$$\{u_1,u_2,u_3\} = \frac{(u_1\bar{u}_2)u_3 + (u_3\bar{u}_2)u_1}{2} \equiv \frac{u_1(\bar{u}_2 u_3) + u_3(\bar{u}_2 u_1)}{2}, \quad (7)$$

$$[u_1,u_2,u_3] = \frac{(u_1\bar{u}_2)u_3 - u_3(\bar{u}_2 u_1)}{2} \equiv \frac{u_1(\bar{u}_2 u_3) - (u_3\bar{u}_2)u_1}{2}, \quad (8)$$

$$\langle u_1, u_2, u_3 \rangle = \frac{(u_1 \bar{u}_2) u_3 - u_1 (\bar{u}_2 u_3)}{2} \equiv \frac{u_3 (\bar{u}_2 u_1) - (u_3 \bar{u}_2) u_1}{2}. \qquad (9)$$

wherein the terms on the right-hand side of (9), besides the alternative order of arranging the parentheses, are also replaced with each other to fit expression (5).

The anticommutator $\{u_1, u_2, u_3\}$ is expressed by a linear combination of arguments, and the commutator $[u_1, u_2, u_3]$ is expressed by a linear combination of the unit $i_0$ and conventional cross products of arguments. So, the definitions (7)-(8) presented as follows:

$$\{u_1, u_2, u_3\} = (u_1, u_2) u_3 - (u_1, u_3) u_2 + (u_2, u_3) u_1, \qquad (10)$$

$$[u_1, u_2, u_3] = ([u_1, u_2], u_3) i_0 - (u_1, i_0)[u_2, u_3] + (u_2, i_0)[u_1, u_3] - (u_3, i_0)[u_1, u_2]. \qquad (11)$$

Formula (10) for $\{u_1, u_2, u_3\}$ is derived from equality (1) by substitution $u_1 + u_3$ in place of $u_1$. Expression (11) for $[u_1, u_2, u_3]$ is derived directly from definition (9) and expression (4) for the product of two hypercomplex numbers.

**5. Interpretation.** Interpretation of the summands $\{u_1, u_2, u_3\}$, $[u_1, u_2, u_3]$ and $\langle u_1, u_2, u_3 \rangle$ in the additive decomposition (5)-(6) of the product $(u_1 \bar{u}_2) u_3$ of three octonions is described by the following.

According to definition (7) and formula (10), $\{u_1, u_2, u_3\}$ is an anticommutator, which is expressed by a linear combination of the arguments $u_1$, $u_2$, $u_3$ and does not change when the first and last arguments, namely $u_1$, and $u_3$, are interchanged. When the multiplicative unit $i_0$ substitutes the central argument, the anticommutator $\{u_1, u_2, u_3\}$ converts into the anticommutator $\{u_1, u_3\}$ for the product of two hypercomplex numbers, described by (3).

The commutator $[u_1, u_2, u_3]$, defined in (8) and represented in (11) as a linear combination of multiplicative unit $i_0$ and conventional cross products of pairs of arguments, is a generalization of the conventional cross product to the case of three arguments both for quaternions and octonions. When one of the arguments is replaced by multiplicative unit $i_0$, the triple cross product $[u_1, u_2, u_3]$ up to a sign coincides with the conventional cross product of two other arguments. When the multiplicative unit $i_0$ substitutes the central argument $u_2$, the commutator $[u_1, u_2, u_3]$ is turned into the commutator $[u_1, u_3]$ for the product of two hypercomplex numbers.

The associator $\langle u_1, u_2, u_3 \rangle$ is well known from the literature [4]. Here it is introduced according to formulas (9), where, unlike conventional definitions, the conjugate central argument is used, and the associator $\langle u_1, u_2, u_3 \rangle$ is computed with a coefficient $\frac{1}{2}$. The mentioned differences do not affect the basic properties of the associator.

**6. Key properties.** As in the case of the cross product of two vectors, the triple cross product $[u_1, u_2, u_3]$ and the associator $\langle u_1, u_2, u_3 \rangle$ have the following common properties:

The triple cross product and the associator are antisymmetric, i.e. when two arguments are the same, cross product vanishes to zero:
$$([u_1, u_2, u_3], u_1) = ([u_1, u_2, u_3], u_2) = ([u_1, u_2, u_3], u_3) = 0,$$
$$(\langle u_1, u_2, u_3 \rangle, u_1) = (\langle u_1, u_2, u_3 \rangle, u_2) = (\langle u_1, u_2, u_3 \rangle, u_3) = 0.$$

Then scalar quadruple mixed product of vectors $([u_1, u_2, u_3], u_4)$ and $(\langle u_1, u_2, u_3 \rangle, u_4)$ are antisymmetric and change sign when any two of the four arguments exchange by places, for example:
$$([u_1, u_2, u_3], u_4) = -([u_4, u_2, u_3], u_1),$$
$$(\langle u_1, u_2, u_3 \rangle, u_4) = -(\langle u_4, u_2, u_3 \rangle, u_1).$$
Here the latter properties are a consequence of the previous two.

In comparison with the triple cross product $[u_1, u_2, u_3]$, the associator $\langle u_1, u_2, u_3 \rangle$ satisfies more orthogonality and zeroing conditions. At that, in addition to orthogonality to each of its arguments, the associator $\langle u_1, u_2, u_3 \rangle$ is orthogonal to the multiplicative unit $i_0$ and to the conventional cross products of its arguments. In total, the orthogonality of the associator to seven vectors is guaranteed:
$$(\langle u_1, u_2, u_3 \rangle, u_1) = (\langle u_1, u_2, u_3 \rangle, u_2) = (\langle u_1, u_2, u_3 \rangle, u_3) = 0,$$
$$(\langle u_1, u_2, u_3 \rangle, i_0) = (\langle u_1, u_2, u_3 \rangle, [u_1, u_2]) = (\langle u_1, u_2, u_3 \rangle, [u_1, u_3]) = (\langle u_1, u_2, u_3 \rangle, [u_2, u_3]) = 0.$$

In general case, the associator $\langle u_1, u_2, u_3 \rangle$ vanishes to zero if its arguments belong to the same quaternion subalgebra, for example:
$$\langle i_0, u_2, u_3 \rangle = \langle [u_2, u_3], u_2, u_3 \rangle = 0.$$

The squares of the lengths of the anticommutator $\{u_1,u_2,u_3\}$, the commutator $[u_1,u_2,u_3]$ and the associator $\langle u_1,u_2,u_3\rangle$ are expressed in terms of the inner products of the vectors $u_1$, $u_2$ and $u_3$ as follows:

$$\|\{u_1,u_2,u_3\}\|^2 = (u_1,u_1)(u_2,u_2)(u_3,u_3) - \det\begin{bmatrix}(u_1,u_1) & (u_1,u_2) & (u_1,u_3)\\(u_1,u_2) & (u_2,u_2) & (u_2,u_3)\\(u_1,u_3) & (u_2,u_3) & (u_3,u_3)\end{bmatrix}, \quad (12)$$

$$\|[u_1,u_2,u_3]\|^2 = ([u_1,u_2],u_3)^2 +$$
$$+\det\begin{bmatrix}(u_1,u_1) & (u_1,u_2) & (u_1,u_3)\\(u_1,u_2) & (u_2,u_2) & (u_2,u_3)\\(u_1,u_3) & (u_2,u_3) & (u_3,u_3)\end{bmatrix} - \det\begin{bmatrix}(u'_1,u'_1) & (u'_1,u'_2) & (u'_1,u'_3)\\(u'_1,u'_2) & (u'_2,u'_2) & (u'_2,u'_3)\\(u'_1,u'_3) & (u'_2,u'_3) & (u'_3,u'_3)\end{bmatrix}, \quad (13)$$

$$\|\langle u_1,u_2,u_3\rangle\|^2 = \begin{vmatrix}(u'_1,u'_1) & (u'_1,u'_2) & (u'_1,u'_3)\\(u'_1,u'_2) & (u'_2,u'_2) & (u'_2,u'_3)\\(u'_1,u'_3) & (u'_2,u'_3) & (u'_3,u'_3)\end{vmatrix} - ([u_1,u_2],u_3)^2. \quad (14)$$

wherein $\|u\|^2$ for any hypercomplex vector $u$ is just a re-designation of $(u,u)$. So that in (12)-(14):

$$\|\{u_1,u_2,u_3\}\|^2 \equiv (\{u_1,u_2,u_3\},\{u_1,u_2,u_3\}),$$
$$\|[u_1,u_2,u_3]\|^2 \equiv ([u_1,u_2,u_3],[u_1,u_2,u_3]),$$
$$\|\langle u_1,u_2,u_3\rangle\|^2 \equiv (\langle u_1,u_2,u_3\rangle,\langle u_1,u_2,u_3\rangle).$$

Using (3), it is easy to verify that in the particular case of the multiplicative unit $i_0$ substituted in (12) as the central factor:

$$\|\{u_1,u_2\}\|^2 = \|\{u_1,i_0,u_2\}\|^2 = (u_1,u_2)^2 - 2(u_1,i_0)(u_2,i_0)(u_1,u_2) +$$
$$+ (u_1,i_0)^2(u_2,u_2) + (u_2,i_0)^2(u_1,u_1)$$

Substituting multiplicative unit $i_0$ in (13) as the central factor (3), we obtain:

$$\|[u_1,u_2]\|^2 = \|[u_1,i_0,u_2]\|^2 = (u_1,u_1)(u_3,u_3) - (u_1,u_3)^2 -$$
$$- (u_1,i_0)^2(u_3,u_3) + 2(u_1,i_0)(u_3,i_0)(u_1,u_3) - (u_3,i_0)^2(u_1,u_1)$$

REMARK: pay attention to the variability of the presentation of the squared lengths, which should be systematized in future works.

Summing (12)-(14) and taking (6) into account, we obtain:

$$\|(u_1\bar{u}_2)u_3\|^2 \equiv ((u_1\bar{u}_2)u_3, (u_1\bar{u}_2)u_3) = (u_1, u_1)(u_2, u_2)(u_3, u_3),$$

matching the axiomatic property of *normalized* algebras for any hypercomplex vectors $u_1, u_2$:

$$(u_1u_2, u_1u_2) = (u_1, u_1)(u_2, u_2).$$

It is characteristic that the triple cross product $[u_1, u_2, u_3]$ changes sign, when conjugated with simultaneous conjugation of arguments, while the anticommutator $\{u_1, u_2, u_3\}$ and the associator $\langle u_1, u_2, u_3 \rangle$ remain unchanged:

$$\overline{\{\bar{u}_1, \bar{u}_2, \bar{u}_3\}} = \{u_1, u_2, u_3\},$$

$$\overline{\langle \bar{u}_1, \bar{u}_2, \bar{u}_3 \rangle} = \langle u_1, u_2, u_3 \rangle,$$

$$\overline{[\bar{u}_1, \bar{u}_2, \bar{u}_3]} = -[u_1, u_2, u_3].$$

Applying the operation of general conjugation with simultaneous conjugation of the arguments to (6), we obtain the decomposition of the product $u_3(\bar{u}_2 u_1)$ of three hypercomplex numbers in the form:

$$u_3(\bar{u}_2 u_1) = \{u_1, u_2, u_3\} - [u_1, u_2, u_3] + \langle u_1, u_2, u_3 \rangle, \tag{15}$$

that was taken into account in the development of the definition (8) for the triple cross product $[u_1, u_2, u_3]$ and in the study of the permutation invariance of mirror symmetry [6,7].

**7. Okubo solution.** In [4] Susumu Okubo considered the product $(u_1 u_2)u_3$ of three octonions $u_1$, $u_2$, $u_3$ and solved the problem of selection of antisymmetric additive term $[u_1, u_2, u_3]_{Ocubo}$. The solution is obtained through cumbersome calculations and is given in [4] on the page 22 by the following expression:

$$(u_1 u_2)u_3 = [u_1, u_2, u_3]_{Ocubo} + 2(u_2, i_0)u_1 u_3 - (u_2, u_3)u_1 - (u_1, u_2)u_3 + (u_1, u_3)u_2,$$

wherein $[u_1, u_2, u_3]_{Okubo}$ presented as follows:

$$[u_1, u_2, u_3]_{Okubo} = \\ = -\langle u_1, u_2, u_3 \rangle + (u_1, i_0)[u_2, u_3] + (u_2, i_0)[u_3, u_1] + (u_3, i_0)[u_1, u_2] - (u_3, [u_1, u_2])i_0.$$

To simplify the comparison of the decomposition of the product of three octonions $(u_1 \bar{u}_2)u_3$ in this paper with the decomposition of three octonions $(u_1 u_2)u_3$ according to [4], the notation matching is established in the Table.

**Table. Notation matching**

| Notion | Notation | |
|---|---|---|
| | in [4] | in this paper |
| Multiplicative unit | $e$ | $i_0$ |
| Inner product | | |
| Conventional cross product of two vectors | $\frac{1}{2}\langle u_1 \vert u_2 \rangle$ | $[u_1, u_2]$ |
| Associator | $-\frac{1}{2}(u_1, u_2, u_3)$ | $\langle u_1, u_2, u_3 \rangle$ |

Taking into account the expressions (10) for $\{u_1, u_2, u_3\}$ and (11) for $[u_1, u_2, u_3]$, the solution [4] is expressed by the formula:

$$(u_1 u_2)u_3 = 2(u_2, i_0)u_1 u_3 - \{u_1, u_2, u_3\} - [u_1, u_2, u_3] - \langle u_1, u_2, u_3 \rangle,$$

which is equivalent to formula (6) for the decomposition of $(u_1 \bar{u}_2)u_3$ into the sum of the anticommutator $\{u_1, u_2, u_3\}$, commutator $[u_1, u_2, u_3]$ and associator $\langle u_1, u_2, u_3 \rangle$.

Thus, formulas (6)-(11) allow us in an obvious way to formulate the idea of decomposing $(u_1 \bar{u}_2)u_3$ into the sum of mutually orthogonal $\{u_1, u_2, u_3\}$, $[u_1, u_2, u_3]$ and $\langle u_1, u_2, u_3 \rangle$ by commutation and changing the order of multiplication of arguments, and formulas (10)-(11) ensures by minimum effort to enhance the results, originally formulated in [4] for $(u_1 u_2)u_3$.

**8. Conclusion.** Unique multiplicative properties of advanced hypercomplex numbers, namely quaternions and octonions, ensure their wide application not only in mathematics and physics, but also in modern information technologies for solving problems of navigation, robotics, computer graphics, image processing, etc.

A well-known method of algebraic computation in terms of hypercomplex numbers is the additive decomposition of the product of two numbers into mutually orthogonal symmetric and antisymmetric components i.e. into anticommutator and a commutator (cross product).

To develop the tools of hypercomplex numbers and their applications, it is interesting to generalize the additive expansion for the product of three arguments. Susumu Okubo proposed an expansion of the product of three octonions [4], wherein the associator is isolated, but the remaining components are not expressed through mutually orthogonal anticommutator and commutator (cross product) generalized to the case of three arguments.

Perhaps, it is just not enough simple form of representation of Okubo's remarkable result [4] for triple product $(u_1 u_2) u_3$ that prevents its widest application in the tools of hypercomplex numbers. It seems possible that the transparent idea of the decomposition of the product $(u_1 \bar{u}_2) u_3$ of three hypercomplex numbers into the sum of mutually orthogonal anticommutator, commutator and associator will help to remove this obstacle and to activate the utilization of octonions.

**Kharinov Mikhail Vyacheslavovich** — PhD., associate professor; senior researcher, Laboratory of Applied Informatics, SPIIRAS. Research interests: digital information analysis, information quantity estimation, numerical representation system, hierarchical data structures, image segmentation and invariant representation for recognition purposes, color transformations of images. The number of publications — 140. khar@iias.spb.su; SPIIRAS, 39, 14-th Line V.O., St. Petersburg, 199178, Russia; office phone +7(812)328-1919, fax +7(812)328-4450.